\documentclass[12pt,reqno]{article}

\usepackage[usenames]{color}
\usepackage{amssymb}
\usepackage{amsmath}
\usepackage{amsthm}
\usepackage{amsfonts}
\usepackage{amscd}
\usepackage{graphicx}

\usepackage[colorlinks=true,
linkcolor=webgreen,
filecolor=webbrown,
citecolor=webgreen]{hyperref}

\definecolor{webgreen}{rgb}{0,.5,0}
\definecolor{webbrown}{rgb}{.6,0,0}

\usepackage{color}
\usepackage{fullpage}
\usepackage{float}

\usepackage{graphics}
\usepackage{latexsym}
\usepackage{epsf}
\usepackage{breakurl}
\usepackage{enumitem}

\newcommand{\seqnum}[1]{\href{https://oeis.org/#1}{\rm \underline{#1}}}

\makeatletter
\def\Ddots{\mathinner{\mkern1mu\raise\p@
		\vbox{\kern7\p@\hbox{.}}\mkern2mu
		\raise4\p@\hbox{.}\mkern2mu\raise7\p@\hbox{.}\mkern1mu}}
\newcommand*\bigcdot{\mathpalette\bigcdot@{.5}}
\newcommand*\bigcdot@[2]{\mathbin{\vcenter{\hbox{\scalebox{#2}{$\m@th#1\bullet$}}}}}
\makeatother

\begin{document}

\theoremstyle{plain}
\newtheorem{theorem}{Theorem}
\newtheorem{corollary}[theorem]{Corollary}
\newtheorem{lemma}[theorem]{Lemma}
\newtheorem{proposition}[theorem]{Proposition}
\newtheorem{claim}[theorem]{Claim}

\theoremstyle{definition}
\newtheorem{definition}[theorem]{Definition}
\newtheorem{example}[theorem]{Example}
\newtheorem{conjecture}[theorem]{Conjecture}

\theoremstyle{remark}
\newtheorem{remark}[theorem]{Remark}

\begin{center}
	\vskip 1cm{\LARGE\bf The Growth Rate of Gijswijt's Sequence
	}
	\vskip 1cm
	\large
	Levi van de Pol\\
	Department of Mathematics\\
	Utrecht University\\
	P.O. Box 80010\\
	3508 TA Utrecht\\
	The Netherlands\\
	\href{mailto:levipol@hotmail.com}{\tt levipol@hotmail.com} \\
\end{center}

\vskip .2 in
\begin{abstract}
Gijswijt's sequence consists almost entirely of small positive integers. However, it is known that every positive integer eventually appears in the sequence. In this paper we determine its growth rate. Specifically, we prove that for $n=4,5,6,\dots$, the number $n$ occurs for the first time at position \mbox{$2\uparrow (2\uparrow(3\uparrow(4\uparrow(5\uparrow\cdots\uparrow((n-2)\uparrow \alpha)))))$}, where $\uparrow$ denotes exponentiation, and $\alpha\in(n-2,n-1)$ is a real number. Our result confirms the growth rate conjectured by van de Bult et al.
\end{abstract}

\section{Introduction}\label{introduction}
To define Gijswijt's sequence, we begin with the concept of a \emph{curling number}, introduced by van de Bult, Gijswijt, Linderman, Sloane, and Wilks~\cite{MR2268453}. For a finite nonempty string $S$ over an arbitrary alphabet, the curling number of $S$ is the largest integer $k$ such that $$S=X\underbrace{YY\cdots Y}_{k\text{ copies}}$$ for some strings $X$ and $Y$, where the string $Y$ must be nonempty. The curling number of $S$ is denoted by $\mathcal{C}(S)$. For example, we have $\mathcal{C}(112112223)=1$. We also have $\mathcal{C}(222322232223)=3$, since the string $2223$ is repeated $3$ times. Using this concept, van de Bult et al.~constructed a recursive sequence $A$:
$$A(1):=1,\qquad 
A(n+1):=\mathcal{C}(A(1),A(2),\dots,A(n))\text{ for all }n\geq1.$$
This sequence is now known as \emph{Gijswijt's sequence}, after Dion Gijswijt, who first defined it~\cite[p.\ 4]{MR2268453}.
The sequence starts with $$A=1,1,2,1,1,2,2,2,3,1,1,2,1,1,2,2,2,3,2,1,1,2,1,\dots $$
Gijswijt's sequence can be found as sequence~\seqnum{A090822} in the On-Line Encyclopedia of Integer Sequences (OEIS)~\cite{oeis}. This sequence consists almost entirely of small positive integers. The first time that a $4$ occurs is at position $220$, and the sequence continues with elements $1$, $2$, $3$, and $4$ for millions of entries. However, van de Bult et al.~\cite[p.\ 12]{MR2268453} proved that all positive integers must occur at some point in Gijswijt's sequence. After this, they used heuristic arguments to specify the approximate location of the first occurrence of a number $n$ in the sequence~\cite[Section 4]{MR2268453}. For $n=1,2,3,4$ this position is respectively $1,3,9,220$. These positions form sequence~\seqnum{A091409} in the OEIS. Surprisingly, the first $5$ occurs at about position $10^{10^{\scriptstyle 23}}$~\cite[p.\ 22]{MR2268453}. The decimal (or binary) representation of this number would require zettabytes of storage. 

The heuristic arguments by van de Bult et al.\ led up to a conjecture~\cite[Conjecture~4.4]{MR2268453}. It states that for $n\geq5$, the first occurrence of $n$ in Gijswijt's sequence is about position
\begin{equation}\label{simple tower}
	2^{2^{\scriptstyle 3^{\scriptstyle 4^{\Ddots^{^{\scriptstyle n-1}}}}}}.
\end{equation}
In this tower of exponents, the 2 occurs twice at the base, followed by the exponents $3,4,\dots,n-2,n-1$. From Expression~\ref{simple tower} we see that Gijswijt's sequence has a very slow growth rate.

The goal of this paper is to establish a precise version of the conjecture. We do this by providing explicit upper and lower bounds for the position of the first $n$ in the sequence. Our main result is as follows:

\begin{theorem}\label{full tower}
For all integers $n\geq4$, the first occurrence of $n$ in {Gijswijt}'s sequence is at position
\begin{equation*}2^{2^{\scriptstyle 3^{\scriptstyle 4^{\Ddots^{^{\scriptstyle (n-2)^{\scriptstyle \alpha}}}}}}},
\end{equation*}
for a real number $\alpha\in (n-2,n-1)$ that depends on $n$.
\end{theorem}

Take, for example, the number $n=4$. The first occurrence in Gijswijt's sequence is at position $220 = 2^{2^{\scriptstyle 2.96\dots}}$, so $\alpha=2.96\ldots\in(2,3)$.
Now consider the number $n=6$. It follows from our theorem that the first $6$ occurs between the positions $2^{2^{\scriptstyle 3^{\scriptstyle 4^{\scriptstyle 4}}}}\approx10^{10^{10^{121}}}$ and $2^{2^{\scriptstyle 3^{\scriptstyle 4^{\scriptstyle 5}}}}\approx10^{10^{\scriptstyle 10^{\scriptstyle 489}}}$.
On the other hand, Theorem~\ref{full tower} is false for $n=3$. Gijswijt's sequence starts with $1,1,2,1,1,2,2,2,3\dots$, so the first $3$ occurs at position $9=2^{3.1699\dots}$. The exponent $3.1699\dots$ falls outside the interval $(1,2)$.

This paper is devoted to the proof of Theorem~\ref{full tower}. We do this as follows. In Section~\ref{the recursive structure} we recall the `recursive structure' that was found by van de Bult et al.~\cite[Section 2]{MR2268453}. This structure makes Gijswijt's sequence $A=A^{(1)}$ the first member of a family of sequences $A^{(m)}$. These sequences are built from finite strings called blocks and glue strings. In Section~\ref{the function t} we show how the positions of first occurrence of $n$ in the various $A^{(m)}$ relate to each other. In Sections~\ref{some inequalities} and~\ref{the asymptotic behavior of beta} we derive inequalities involving the lengths of the blocks and glue strings. In Section~\ref{the tower of exponents} we prove our main result by combining the inequalities into a single expression. Remark~\ref{narrow down} explains a strategy to obtain sharper bounds.

Our Section~\ref{the function t} is based on the ideas in Section 4.6 by van de Bult et al.~\cite{MR2268453}. At first glance, the arguments in our Section~\ref{some inequalities} have nothing to do with their heuristic derivation~\cite[Section 4]{MR2268453}. However, our Section~\ref{some inequalities} is actually inspired by the `tabular construction' and `smoothing' techniques that they provided. In their Section~4, van de Bult et al.~made three other conjectures about related inequalities and the recursive structure~\cite[Conjectures 4.1-4.3]{MR2268453}. These conjectures served mainly as stepping stones towards Expression~\ref{simple tower}. It might be possible to derive the conjectures from our work, but we prefer to focus on the main result.

\subsection{The curling number conjecture}\label{cnc}
Van de Bult et al.~\cite{MR2268453} posed another conjecture near the end of their paper~\cite[Finiteness Conjecture]{MR2268453}. In a different paper, Chaffin et al.~\cite{MR3056627} called this conjecture the \emph{curling number conjecture}.
\begin{conjecture}[The curling number conjecture~{\cite[Conjecture 1]{MR3056627}}]
	If one starts with any initial string of integers
	$S$, and extends it by repeatedly appending the curling number of the current string, the string
	will eventually reach 1.
\end{conjecture}
The techniques in this paper are completely unsuitable for proving this conjecture. The curling number conjecture seems harder than the one regarding the growth rate of Gijswijt's sequence. However, Chaffin et al.~\cite{MR3056627} made some progress. Among other things, they computationally verified the conjecture for all initial strings with lengths up to 48 that consist only of the digits $2$ and $3$.

\subsection{Notation}
We use the term \emph{string} for a finite list $e_1,e_2,\dots,e_l$ of positive integers and reserve
the term \emph{sequence} for infinite lists $e_1, e_2,\dots$ of positive integers. The objects $e_1,e_2,\dots,e_l$ are the \emph{elements} of the string. The \emph{length} of a string is the number of elements, including duplicates. The empty string $\epsilon$ has length $0$. We usually omit the commas, so 1222 means $1,2,2,2$. We can \emph{concatenate} strings $S=s_1\cdots s_l$ and $T=t_1\cdots t_m$ to obtain a new string $W=ST=s_1\cdots s_lt_1\cdots t_m$. Similarly, the concatenation $W = ST$ is a sequence when $S$ is a string and $T$ is a sequence. In either case, the string $S$ is a \emph{prefix} of $W$ and a \emph{proper prefix} if $T$ is not empty. When $W = ST$ is a string, the string $T$ is a \emph{suffix} of $W$ and a \emph{proper suffix} if $S$ is not empty. When $W = USV$ is a string or sequence, the string $S$ is a \emph{substring} of $W$. We denote the $k$th element of $S$ by $S(k)$. Consider strings $S_1,S_2,S_3,\dots$ such that $S_n$ is a proper prefix of $S_{n+1}$ for all $n\geq1$. We define the \emph{limit} $\lim_{n\rightarrow\infty}S_n$ as the unique sequence $S$ such that $S_n$ is a prefix of $S$ for all $n\geq1$.

For a string $X$ and an integer $k\geq1$, we denote $\underbrace{XX\cdots X}_{k\text{ copies}}$ by $X^k$. Now the \emph{curling number} of a string $S$ is the largest integer $c$ such that $S$ can be written as $XY^c$ with $Y$ nonempty.
Van de Bult et al.~\cite{MR2268453} only defined the curling number of a nonempty string; we additionally define $\mathcal{C}(\epsilon):=1$. 
Throughout this paper, the expression $(m+1)^k$ can both denote the string $(m+1)^k=\underbrace{(m+1)(m+1)\cdots (m+1)}_{k\text{ copies}}$ or the number $(m+1)^k=\underbrace{(m+1)\cdot(m+1)\cdot\ldots\cdot (m+1)}_{k\text{ copies}}$. The difference will be clear from the context. We always use $\cdot$ for multiplication, and no symbol for concatenation.

\section{The recursive structure}\label{the recursive structure}
Van de Bult et al.\ defined a family of sequences that generalize Gijswijt's sequence~\cite[Section~3]{MR2268453}. For $m\geq1$, the $m$th order sequence $A^{(m)}$ is given by
$$A^{(m)}(1):=m,\qquad 
A^{(m)}(n+1):=\max(m,\mathcal{C}(A^{(m)}(1),A^{(m)}(2),\dots,A^{(m)}(n)))\text{ for all }n\geq1.$$
Note that $A^{(1)}=A$.
The sequences $A^{(2)}$ and $A^{(3)}$ 
start as follows:
$$A^{(2)}=2,2,2,3,2,2,2,3,2,2,2,3,3,2,2,2,3,2,2,2,3,2,\dots$$
$$A^{(3)}=3,3,3,3,4,3,3,3,3,4,3,3,3,3,4,3,3,3,3,4,4,3,\dots$$
These are sequences~\seqnum{A091787} and~\seqnum{A091799} in the OEIS~\cite{oeis}. Van de Bult et al.\ proved that each sequence $A^{(m)}$, for $m\geq1$, is built from finite strings $B_t^{(m)}$ and $S_t^{(m)}$. The strings $B_t^{(m)}$ are called \emph{blocks} and the strings $S_t^{(m)}$ are called \emph{glue strings}. 
They satisfy the following relations for $t,m\geq1$:
\begin{align}
	B_1^{(m)}&=m\nonumber\\
B_{t+1}^{(m)}&=\underbrace{B_t^{(m)}B_t^{(m)}\cdots B_t^{(m)}}_{m+1\text{ copies}}S_t^{(m)}	\nonumber\\
A^{(m)}&=\lim_{t\rightarrow\infty}B_t^{(m)}\nonumber\\
A^{(m+1)}&=S_1^{(m)}S_2^{(m)}S_3^{(m)}\cdots\nonumber
\end{align}
The limit in the third line means that every block $B_t^{(m)}$ is a prefix of $A^{(m)}$. Therefore we can obtain the sequence $A^{(m)}$ by starting with $B_1^{(m)}$, and then repeatedly applying the second equation. In each iteration, we copy the current block $B_t^{(m)}$ for $m+1$ times, and then add the glue string $S_t^{(m)}$ to find $B_{t+1}^{(m)}$. Surprisingly, the glue strings used for constructing $A^{(m)}$ can be combined to form the sequence $A^{(m+1)}$; this is the fourth equation. Note that the blocks $B_t^{(m)}$ contain only elements $m$ and higher, while the glue strings $S_t^{(m)}$ contain only elements $m+1$ and higher.

A few examples are in order. Consider the sequence $A^{(2)}$. The first prefix blocks are $B_1^{(2)}=2$, $B_2^{(2)}=2223$, $B_3^{(2)}=(2223)^33$, $B_4^{(2)}=((2223)^33)^3334$. The associated glue strings are $S_1^{(2)}=3$, $S_2^{(2)}=3$, $S_3^{(2)}=334$. When we concatenate these glue strings we get $33334$, which is the prefix $B_2^{(3)}$ of $A^{(3)}$.
In turn, the block $B_2^{(3)}$ equals $(B_1^{(3)})^4S_1^{(3)}$. We have $B_1^{(3)}=3$ and $S_1^{(3)}=4$. The string $S_1^{(3)}=4$ is the prefix $B_1^{(4)}$ of $A^{(4)}$.

This structure implies that every element of $A^{(m+1)}$ also occurs in $A^{(m)}$. Indeed, if $x$ is in $A^{(m+1)}$, then there is a $t$ such that $x$ is in $S_t^{(m)}$. Therefore, the element $x$ occurs in $B_{t+1}^{(m)}$, and hence in $A^{(m)}$. Now from induction we see that for $n\geq m$, every element of $A^{(n)}$ occurs in $A^{(m)}$. Since $A^{(n)}$ starts with $n$, we conclude that $A^{(m)}$ contains all integers that are at least $m$. Like Gijswijt's sequence, the sequences $A^{(m)}$ for $m\geq2$ are also slow-growing; see Remark~\ref{higher order}.

Additionally, van de Bult et al.~\cite{MR2268453} defined strings $T_t^{(m)}$ for $t,m\geq1$:
\begin{equation*}T_t^{(m)}:=S_1^{(m)}S_2^{(m)}\cdots S_{t-1}^{(m)}.
\end{equation*}
Now it follows that
$$A^{(m+1)}=\lim_{t\rightarrow\infty}T_t^{(m)}.$$
Since $T_t^{(m)}$ is a prefix of $A^{(m+1)}$ for all $m,t\geq1$, the strings $T_t^{(m)}$ consist of elements $m+1$ and higher. Here is a short table for reference:
\begin{center}
	\begin{table}[H]\centering
		\begin{tabular}{ |p{2.7cm}|p{2.7cm}|p{2.7cm}|p{2.7cm}|p{2.7cm}|}
			\hline
			$B_1^{(m)}$ & $B_2^{(m)}$& $S_1^{(m)}$  & $T_1^{(m)}$ & $T_2^{(m)}$\\
			\hline
			$m$   &$ m^{m+1}(m+1)    $&$m+1$&$ \epsilon$&$m+1$\\
			\hline
		\end{tabular}
		\caption{Blocks and glue strings for small $t$ and general $m$.}
	\end{table}
\end{center}
The strings $B_t^{(m)}$ and $S_t^{(m)}$ are always nonempty. However, the strings $T_1^{(m)}$ are empty for all $m$.
The lengths of $$B_t^{(m)},S_t^{(m)},T_t^{(m)}$$ are denoted by respectively $$\beta^{(m)}(t),\sigma^{(m)}(t),\tau^{(m)}(t).$$ 
For all $m$ and $t$, the string $B_t^{(m)}$ is a proper prefix of $B_{t+1}^{(m)}$ and the string $T_t^{(m)}$ is a proper prefix of $T_{t+1}^{(m)}$. Therefore the functions $\beta^{(m)}(t)$ and $\tau^{(m)}(t)$ are strictly increasing in $t$. This is not true for $\sigma^{(m)}(t)$.

We often use the following property of curling numbers: 
\begin{lemma}\label{basic}
		If $S$ is a string with suffix $T$, then $\mathcal{C}(S)\geq\mathcal{C}(T)$.
	\begin{proof}
		Let $c:= \mathcal{C}(T)$. Then $T$ has a suffix $Y^c$ with $Y$ nonempty. It follows that $S$ has suffix $Y^c$ as well, so $\mathcal{C}(S)\geq c = \mathcal{C}(T)$.
	\end{proof}
\end{lemma}

We can now start our derivation of Theorem~\ref{full tower}. As a first step we prove three lemmas about the recursive structure. 
\begin{lemma}\label{T suffix} 
	For all $m\geq1$ and $t\geq2$, the string $mT_t^{(m)}$ is a suffix of $B_t^{(m)}$. This is the only occurrence of $T_t^{(m)}$ as a substring of $B_t^{(m)}$.
	\begin{proof}
		We prove the first statement by induction on $t$. For $t=2$, the string $mT_2^{(m)}=m(m+1)$ is a suffix of $B_2^{(m)}=m^{m+1}(m+1)$.
		
		Now suppose that the first statement is true for $t=k$. Then we use that $$B_{k+1}^{(m)}=(B_k^{(m)})^{m+1}S_k^{(m)}.$$ By the induction hypothesis, the string $B_{k+1}^{(m)}$ has suffix $mT_k^{(m)}S_k^{(m)}$. The suffix equals $mT_{k+1}^{(m)}$ by definition of $T_{k+1}^{(m)}$. This completes the induction step.
		
		The second statement is Lemma 3.2b by van de Bult et al.~\cite{MR2268453}.
	\end{proof}
\end{lemma}
\begin{lemma}\label{t strings property}
	For all $m\geq1$, the strings $T_t^{(m)}$ are exactly the prefixes of $A^{(m+1)}$ with curling number at most $m$.
	\begin{proof}
		The string $T_1^{(m)}=\epsilon$ is a prefix of $A^{(m+1)}$ with curling number $1$.
		
		For $t\geq2$, the block $B_{t+1}^{(m)}$ is a prefix of $A^{(m)}$. This block equals $(B_t^{(m)})^{m+1}S_t^{(m)}$, so since $B_t^{(m)}(1)=m$, the string $B_t^{(m)}m$ is also a prefix of $A^{(m)}$. Now from the definition of $A^{(m)}$ we obtain that $\max(m,\mathcal{C}(B_t^{(m)}))=m$, so $\mathcal{C}(B_t^{(m)})\leq m$. We know by Lemma~\ref{T suffix} that $T_t^{(m)}$ is a suffix of $B_t^{(m)}$, so Lemma~\ref{basic} tells us that $\mathcal{C}(T_t^{(m)})\leq \mathcal{C}(B_t^{(m)})\leq m$.
		
		Conversely, let $P$ be a nonempty prefix of $A^{(m+1)}$ with curling number at most $m$. Then we can write $P=T_t^{(m)}Q$ for some $t\geq2$ and a proper prefix $Q$ of $S_t^{(m)}$. We have $S_t^{(m)}=QR$ for a nonempty string $R$. All elements of $R$ are at least $m+1$, in particular the first element $R(1)$. Now the string $$(B_t^{(m)})^{m+1}QR(1)$$ is a prefix of $B_{t+1}^{(m)}$, so it is also a prefix of $A^{(m)}$. Therefore $\max(m,\mathcal{C}((B_t^{(m)})^{m+1}Q))=R(1)\geq m+1$, so $\mathcal{C}((B_t^{(m)})^{m+1}Q)\geq m+1$. Hence the string $(B_t^{(m)})^{m+1}Q$ has a nonempty suffix $X^{m+1}$. By the first part of Lemma~\ref{T suffix}, the string $(B_t^{(m)})^{m+1}Q$ also has suffix $mT_t^{(m)}Q$. Since $\mathcal{C}(T_t^{(m)}Q)=\mathcal{C}(P)\leq m$, the string $mT_t^{(m)}Q$ must be a suffix of $X^{m+1}$. Therefore $m\in X^{m+1}$, so $m\in X$. Now it follows that $mT_t^{(m)}Q$ is a suffix of $X$. Therefore $T_t^{(m)}Q$ occurs as a substring of $(B_t^{(m)})^{m+1}$. But $T_t^{(m)}Q$ contains only elements that are at least $m+1$, while $B_t^{(m)}(1)=m$. So $T_t^{(m)}Q$ is a substring of $B_t^{(m)}$. It follows by the second part of Lemma~\ref{T suffix} that the substring $T_t^{(m)}$ must occur at the end of $B_t^{(m)}$, so $Q$ is empty. We conclude that $P$ is of the form $T_t^{(m)}$ with $t\geq2$.	
	\end{proof}
\end{lemma}
\begin{lemma}\label{groups}
	For all $m\geq1$ and $t\geq2$, the elements $m$ in $B_t^{(m)}$ occur in groups of $m+1$ consecutive times $m$. These groups are always followed by an element $m+1$.
\begin{proof} We use induction on $t$. For $t=2$ we have $B_2^{(m)}=m^{m+1}(m+1)$. This is exactly one group of $m+1$ consecutive $m$'s, followed by exactly one element $m+1$. Now suppose that the statement holds for $t=k$. We know that $$B_{k+1}^{(m)}=(B_{k}^{(m)})^{m+1}S_k^{(m)}.$$ Each copy of $B_k^{(m)}$ contains $m$'s in groups of $m+1$ times $m$ followed by an element $m+1$ by the induction hypothesis. The string $S_k^{(m)}$ contains no elements $m$. This finishes the induction step.
\end{proof}
\end{lemma}

\section{The function \texorpdfstring{$t(m,n)$}{t(m,n)}}\label{the function t}
We now turn our attention to the first occurrence of $n$ in the sequence $A^{(m)}$, where $n\geq m$. We show in Lemma~\ref{t existence} that the first $n$ in $A^{(m)}$ occurs at the end of a block $B_t^{(m)}$. The question is then: what is the value of $t$? We denote this value by $t(m,n)$. 
\begin{lemma}\label{t existence}
	For all integers $m,n$ with $n\geq m\geq1$, the first occurrence of $n$ in $A^{(m)}$ is as the last element of a prefix block $B_t^{(m)}$ for some $t\geq1$.
	\begin{proof}
		For $n=m$, the number $m$ occurs at the beginning of $A^{(m)}$. This is also at the end of $B_1^{(m)}=m$. Now suppose that $n>m$. Let $t$ be minimal such that $n$ is in $B_t^{(m)}$. Since $n$ is not an element of $B_{t-1}^{(m)}$, the first $n$ in $A^{(m)}$ must be in the glue string $S_{t-1}^{(m)}$. Let $B$ be the prefix of $A^{(m)}$, ending at the first $n$. Now $(B_{t-1}^{(m)})^{m+1}$ is a proper prefix of $B$, and $B$ is a prefix of $B_t^{(m)}$. Since $n$ appears at the end of $B$ and nowhere else in $B$, we have $\mathcal{C}(B)=1$. Now by definition of $A^{(m)}$, the prefix $B$ of $A^{(m)}$ is followed by $m$. Since $m\not\in S_{t-1}^{(m)}$, it follows that $B=B_t^{(m)}$. We conclude that the first $n$ in $A^{(m)}$ is the last element of $B_t^{(m)}$. 
	\end{proof}
\end{lemma}
Hence we can make the following definition:
\begin{definition}
	For all integers $m,n$ with $n\geq m\geq 1$, define $t(m,n)$ as the integer such that the first $n$ in $A^{(m)}$ occurs at the end of prefix block $B^{(m)}_{t(m,n)}$.
\end{definition}

Van de Bult et al.~\cite[p.\ 22]{MR2268453} deduced that the first occurrence of $5$ in $A^{(2)}$ is at the end of the block $B_{80}^{(2)}$. Therefore it occurs at position $\beta^{(2)}(80)$. In our new terminology we obtain $t(2,5)=80$. In general, the first occurrence of $n$ in $A^{(m)}$ is at position $\beta^{(m)}(t(m,n))$. The values of $t(1,n)$ are given in sequence~\seqnum{A357064} of the OEIS~\cite{oeis}.

The function $t(m,n)$ is not only related to the blocks $B_t^{(m)}$, but also to the strings $T_t^{(m)}$:
\begin{lemma}\label{tm strings equality} For $n\geq m+1$ we have the relation $T^{(m)}_{t(m,n)}=B^{(m+1)}_{t(m+1,n)}$.
	\begin{proof}	
		By Lemma~\ref{T suffix}, the string $T^{(m)}_{t(m,n)}$ is a suffix of $B^{(m)}_{t(m,n)}$. It follows by definition of $t(m,n)$ that the last element of $T^{(m)}_{t(m,n)}$ is $n$, and that $T^{(m)}_{t(m,n)}$ has no other elements $n$. Also, the string $T^{(m)}_{t(m,n)}$ is a prefix of $A^{(m+1)}$. Therefore the first occurrence of $n$ in $A^{(m+1)}$ is both at the end of $T^{(m)}_{t(m,n)}$ and at the end of $B^{(m+1)}_{t(m+1,n)}$; hence the two strings are equal.	
	\end{proof}
\end{lemma}

Let us consider a few examples.
\begin{example}\label{t example}
	For all $m\geq1$, we have the following:
	\begin{itemize} 
		\item $t(m,m)=1$ since $B_1^{(m)}=m$.
		\item  $t(m,m+1)=2$ since $B_2^{(m)}=m^{m+1}(m+1)$.
		\item For $n=m+2$, we need the values $\sigma^{(m)}(x)=1$ for $x=1,2,\dots,m$ and $\sigma^{(m)}(m+1)=3$ by Equation 35 from van de Bult et al.~\cite{MR2268453}. Summing these values, we obtain $\tau^{(m)}(m+2)=m+3$. Now $T^{(m)}_{m+2}$ consists of the first $m+3$ elements of $A^{(m+1)}$. It follows that $T^{(m)}_{m+2}=(m+1)^{m+2}(m+2)=B_2^{(m+1)}=B_{t(m+1,m+2)}^{(m+1)}$. It follows from Lemma~\ref{tm strings equality} that $t(m,m+2)=m+2$.
	\end{itemize}
\end{example} 
\noindent Lemma~\ref{tm strings equality} tells us that $T^{(m)}_{t(m,n)}=B^{(m+1)}_{t(m+1,n)}$. Therefore their lengths are equal:
\begin{equation}\label{tau tau eq}\tau^{(m)}(t(m,n))=\beta^{(m+1)}(t(m+1,n)).\end{equation}
In Section~\ref{some inequalities} we prove inequalities between values of the functions $\tau^{(m)}$ and $\beta^{(m+1)}$. From this we deduce bounds for $t(m,n)$ in terms of $t(m+1,n)$. We can then combine this for various values of $m$, to obtain bounds for $t(1,n)$ in terms of $t(n-2,n)$, which equals $n$ by Example~\ref{t example}. Since the first occurrence of $n$ in Gijswijt's sequence is at position $\beta^{(1)}(t(1,n))$, it remains to bound $\beta^{(1)}(t)$ in terms of $t$. We do that in Section~\ref{the asymptotic behavior of beta}.

\section{Linking \texorpdfstring{$\beta^{(m+1)}$}{beta} and \texorpdfstring{$\tau^{(m)}$}{tau}}\label{some inequalities}
In this section we bound the value of $\beta^{(m+1)}(u)$ by $\tau^{(m)}(t)$ for certain values of $t$ and $u$. We do this with a series of lemmas that leads up to the inequalities in Lemma~\ref{estimates}. Using these inequalities we bound $t(m,n)$ in terms of $t(m+1,n)$. The result can be found in Proposition~\ref{estimates 3}.

Our strategy is as follows. Fix $u,m\geq1$. We know that every string $T^{(m)}_t$ is a prefix of $A^{(m+1)}$. The string $B^{(m+1)}_u$ is also a prefix of $A^{(m+1)}$. We prove lower and upper bounds for the number of proper prefixes of $B^{(m+1)}_u$ that are of the form $T^{(m)}_t$. If there are at least $L$ such proper prefixes, then it follows that $\beta^{(m+1)}(u)\geq \tau^{(m)}(L)$. If there are at most $U$ such proper prefixes, then we get $\beta^{(m+1)}(u)\leq \tau^{(m)}(U+1)$. In Lemma~\ref{estimates} we prove this with $L=m\cdot (m+1)^{u-2}$ and $U=(m+1)^{u-1}$.

In Lemmas~\ref{curl at most m}--\ref{a1 not zero} we analyze the curling numbers of proper prefixes of $B^{(m+1)}_u$, as this determines whether a given prefix can be written as $T_t^{(m)}$ via Lemma~\ref{t strings property}. 
We combine these results in Lemma~\ref{Ttm characterization}. This lemma proves that certain proper prefixes of $B_u^{(m+1)}$ are of the form $T^{(m)}_t$, and gives restrictions on which proper prefixes can be of the form $T^{(m)}_t$. This gives us lower and upper bounds for the amount of proper prefixes of the form $T_t^{(m)}$, which is exactly what we need.

Recall from Lemma~\ref{t strings property} that the strings $T_t^{(m)}$ are exactly the prefixes of $A^{(m+1)}$, and hence of some $B_u^{(m+1)}$, that have curling number at most $m$. Before investigating these strings further, we study the prefixes of $B_u^{(m+1)}$ with curling number at most $m+1$. The following lemma gives a characterization of those strings. 
\begin{lemma}\label{curl at most m}
	For $m,u\geq1$, the proper prefixes $P$ of $B_u^{(m+1)}$ with $\mathcal{C}(P)\leq m+1$ are exactly the strings of the form \begin{equation}(B_{u-1}^{(m+1)})^{a_{u-1}}\cdots(B_2^{(m+1)})^{a_2}(B_1^{(m+1)})^{a_1}\label{eq}
	\end{equation} with $a_1,a_2,\dots,a_{u-1}\in\{0,1,\dots,m+1\}$.
	\begin{proof}
		By definition of $A^{(m+1)}$, the prefixes with curling number at most $m+1$ are the prefixes that are followed by an element $m+1$. Every string $B_u^{(m+1)}$ is a prefix of $A^{(m+1)}$, so the same is true for all proper prefixes of $B_u^{(m+1)}$. We prove by induction on $u$ that these strings have the desired form. For $u=1$, the string $B_1^{(m+1)}=(m+1)$ has only one proper prefix $\epsilon$, which is followed by $m+1$. This corresponds to Expression~\ref{eq} since there are no numbers between $1$ and $u-1$ for $u=1$. 
		
		Now suppose that the statement holds for $u=k$. Then we use that $$B_{k+1}^{(m+1)}=(B_k^{(m+1)})^{m+2}S_k^{(m+1)}.$$ The string $B_k^{(m+1)}$ starts with $m+1$ while $S_k^{(m+1)}$ contains no elements $m+1$. Therefore the proper prefixes $P$ of $B_{k+1}^{(m+1)}$ that are followed by $m+1$ are of the form $$P=(B_k^{(m+1)})^aP',$$ where $0\leq a\leq m+1$ and $P'$ is a proper prefix of $B_k^{(m+1)}$ that is followed by an element $m+1$. By the induction hypothesis, the string $P'$ is of the form $$(B_{k-1}^{(m+1)})^{a_{k-1}}\cdots(B_2^{(m+1)})^{a_2}(B_1^{(m+1)})^{a_1}.$$ Concatenating this to $(B_k^{(m+1)})^a$ completes the induction step.
	\end{proof}
\end{lemma}

To find whether a string $(B_{u-1}^{(m+1)})^{a_{u-1}}\cdots (B_1^{(m+1)})^{a_1}$ can be written as $T_t^{(m)}$ for some $t$, we need to know whether its curling number is at most $m$, or equal to $m+1$. For a large class of strings, we can easily prove that this curling number is exactly $m+1$:
\begin{lemma}\label{maxdigit}
	Let $m,u\geq1$, and consider a prefix $P:=(B_{u-1}^{(m+1)})^{a_{u-1}}\cdots(B_1^{(m+1)})^{a_1}$ of $B_u^{(m+1)}$, for
	$a_1,\dots,a_{u-1}\in\{0,1,\dots,m+1\}$. Suppose that 
	$a_i=m+1$ for some $i$. Then $\mathcal{C}(P)= m+1$.
	\begin{proof}
		By Lemma~\ref{curl at most m}, the string $Q:=(B_{i-1}^{(m+1)})^{a_{i-1}}\cdots(B_1^{(m+1)})^{a_1}$ is a prefix of $B_i^{(m+1)}$. Therefore we can write $B_i^{(m+1)}=QR$. The string $(B_i^{(m+1)})^{a_i} Q$ is a suffix of $P$, and $(B_i^{(m+1)})^{a_i} Q$ is equal to $(QR)^{a_i} Q=Q(RQ)^{a_i}=Q (RQ)^{m+1}$. Therefore, we have $\mathcal{C}(P)\geq m+1$.
			
		By Lemma~\ref{curl at most m}, this curling number is also \emph{at most} $m+1$. Therefore the curling number is exactly $m+1$. 
	\end{proof}
\end{lemma}

It follows from Lemmas~\ref{t strings property} and~\ref{maxdigit} that all strings $(B_{u-1}^{(m+1)})^{a_{u-1}}\cdots(B_1^{(m+1)})^{a_1}$ where at least one of $a_1,\dots,a_{u-1}$ is equal to $m+1$, are not of the form $T^{(m)}_t$. 
The remaining strings that we have to consider are the strings $(B_{u-1}^{(m+1)})^{a_{u-1}}\cdots(B_1^{(m+1)})^{a_1}$ with $a_i\in\{0,1,\dots,m\}$ for all $i$.  
In the following lemma, we establish that under the condition $a_1\neq0$, the curling number of such strings is at most $m$.
\begin{lemma}\label{a1 not zero}
	For $m\geq1$ and $u\geq2$, let $$P= (B_{u-1}^{(m+1)})^{a_{u-1}}\cdots(B_1^{(m+1)})^{a_1}$$ be a prefix of $B_u^{(m+1)}$ with $a_1,\dots,a_{u-1}\in\{0,1,\dots,m\}$ and $a_1\neq0$. Then we have $\mathcal{C}(P)\leq m.$
	\begin{proof}		
		We prove this by contradiction, so suppose that $P$ has a nonempty suffix of the form $X^{m+1}$. We can then write $$P=(B_{u-1}^{(m+1)})^{a_{u-1}}\cdots(B_1^{(m+1)})^{a_1}=QX^{m+1}.$$ Since $a_1\neq 0$, the last element of $P$ is $m+1$. Therefore $m+1$ is also the last element of $X$.
		
		We now prove a claim about the first elements of $X$. 
		\begin{claim}\label{claim 1}
			The string $X$ has prefix $(m+1)^{k}(m+2)$ for some $k>0$.
		\begin{proof}
		Since $a_1\leq m$, the curling number of $(B_1^{(m+1)})^{a_1}=(m+1)^{a_1}$ is smaller than $m+1$. We see that $X^{m+1}$ has to be longer than $(m+1)^{a_1}$ (that is, it has more elements). It follows from Lemma~\ref{T suffix} that all strings $B_i^{(m+1)}$ with $i\geq2$ end with an element larger than $m+1$. Hence $X^{m+1}$ has a suffix $l(m+1)^{a_1}$ for some $l>m+1$. This can only happen when $l(m+1)^{a_1}$ is a suffix of $X$ itself. Now $P$ looks like this:
		
$$P=(B_{u-1}^{(m+1)})^{a_{u-1}}\cdots(B_1^{(m+1)})^{a_1}=Q\underbrace{\ldots l(m+1)^{a_1}}_{\text{first }X}\underbrace{\ldots l(m+1)^{a_1}}_{\text{second }X}\cdots \underbrace{\ldots l(m+1)^{a_1}}_{\text{last }X}$$

By Lemma~\ref{groups}, applied to $B_u^{(m+1)}$, the elements $m+1$ in $P$ occur in groups of $m+2$ times $m+1$, followed by an element $m+2$. Since we cut off the string $B_{u}^{(m+1)}$ at the end of $P$, there is a possible exception to this at the end of $P$. When we now consider that the first copy of $X$ in $X^{m+1}$ has suffix $l(m+1)^{a_1}$, we deduce that the second copy of $X$ (and hence all of them) has prefix $(m+1)^{m+2-a_1}(m+2)$. From $a_1\leq m$ it follows that $k:=m+2-a_1>0$.
		\end{proof}
		\end{claim}
		Next we prove that the first element of $X$ is larger than $m+1$, which contradicts the previous claim.
		\begin{claim}\label{claim 2}
			The first element of $X$ is larger than $m+1$.
		\begin{proof}				
		Define $i$ as the largest integer such that $X^{m+1}$ contains $T_i^{(m+1)}$ as a substring. We saw in Claim~\ref{claim 1} that $X^{m+1}$ contains an element $m+2=T_2^{(m+1)}$. Therefore $i\geq2$. Since $X$ ends with $m+1$ and $T_i^{(m+1)}$ contains no element $m+1$, the string $T_i^{(m+1)}$ has to be a substring of $X$. Therefore $X^{m+1}$ contains at least $m+1$ disjoint copies of $T_i^{(m+1)}$. 
		
		By Lemma~\ref{curl at most m}, the string $P_{i-1}:=(B_{i-1}^{(m+1)})^{a_{i-1}}\cdots(B_1^{(m+1)})^{a_1}$ is a proper prefix of $B_i^{(m+1)}$. It follows from Lemma~\ref{T suffix} that $P_{i-1}$ does not have $T_i^{(m+1)}$ as a substring. Again using Lemma~\ref{T suffix}, the string $P_i:=(B_i^{(m+1)})^{a_i}\cdots(B_1^{(m+1)})^{a_1}$ contains exactly $a_i$ copies of $T_i^{(m+1)}$. Now since $a_i\leq m$, the string $X^{m+1}$ contains more copies of $T_i^{(m+1)}$ than $P_i$. Both $X^{m+1}$ and $P_i$ are suffixes of $P$. Therefore $P_i$ is a proper suffix of $X^{m+1}$.
		
		Let $j$ be the smallest integer larger than $i$ such that $a_j\neq0$. Then $B_j^{(m+1)}P_i$ is a suffix of $P$. The string $B_j^{(m+1)}$ has suffix $T_j^{(m+1)}$ by Lemma~\ref{T suffix}, but $X^{m+1}$ cannot contain $T_j^{(m+1)}$ as a substring by definition of $i$. So $X^{m+1}$ is a proper suffix of $T_j^{(m+1)}P_i$. Remembering that $P_i$ is a proper suffix of $X^{m+1}$, it follows that the first element of $X^{m+1}$ lies inside $T_j^{(m+1)}$. Hence this element is larger than $m+1$.
		\end{proof}
		\end{claim}
\noindent	
Since the two claims contradict each other, the lemma must be true.
	\end{proof} 
\end{lemma}

Combining the above lemmas, we conclude the following about the strings $T_t^{(m)}$:
\begin{lemma}\label{Ttm characterization}
	Let $m\geq1$ and $u\geq2$. 
	\begin{enumerate}[label=(\alph*)] 
		\item\label{a} For each proper prefix $T_t^{(m)}$ of $B^{(m+1)}_u$, there are $a_1,\dots,a_{u-1}\in\{0,1,\dots,m\}$ such that $T_t^{(m)}=(B_{u-1}^{(m+1)})^{a_{u-1}}\cdots(B_1^{(m+1)})^{a_1}$. 
		\item\label{b} For all integers $a_1,\dots,a_{u-1}\in\{0,1,\dots,m\}$ with $a_1\neq0$, there is a proper prefix $T_t^{(m)}$ of $B^{(m+1)}_u$ such that $(B_{u-1}^{(m+1)})^{a_{u-1}}\cdots(B_1^{(m+1)})^{a_1}=T_t^{(m)}$.
	\end{enumerate}		
	\begin{proof}\hfill
		\begin{enumerate} [label=(\alph*)] 
			\item By Lemma~\ref{t strings property}, the string $T_t^{(m)}$ has curling number at most $m$. Now it follows from Lemma~\ref{curl at most m} that $T_t^{(m)}$ can be written as $(B_{u-1}^{(m+1)})^{a_{u-1}}\cdots(B_1^{(m+1)})^{a_1}$ with $a_1,\dots,a_{u-1}\in\{0,\dots,m+1\}$. By Lemma~\ref{maxdigit}, all these exponents are at most $m$.
			\item Consider a string $(B_{u-1}^{(m+1)})^{a_{u-1}}\cdots(B_1^{(m+1)})^{a_1}$ with $m\geq1$, $a_1\neq0$, and $a_1,\dots,a_{u-1}\in\{0,\dots,m\}$. By Lemma~\ref{curl at most m}, this string is a proper prefix of $B_u^{(m+1)}$, and hence of $A^{(m+1)}$. By Lemma~\ref{a1 not zero}, the string has curling number at most $m$. Now it follows from Lemma~\ref{t strings property} that $(B_{u-1}^{(m+1)})^{a_{u-1}}\cdots(B_1^{(m+1)})^{a_1}$ is equal to $T_t^{(m)}$ for some $t$.
		\end{enumerate}
	\end{proof}
\end{lemma}

We can now deduce inequalities between values of $\beta^{(m+1)}$ and $\tau^{(m)}$. All further estimates are based on the following lemma.
\begin{lemma}\label{estimates}
	For all $m\geq1$ and $u\geq2$ we have $$\tau^{(m)}(m\cdot (m+1)^{u-2})\leq \beta^{(m+1)}(u)\leq \tau^{(m)}((m+1)^{u-1}+1).$$ When $u=m+3$, the second inequality is strict.
\begin{proof}
	By Lemma~\ref{Ttm characterization}~\ref{b}, the string $B^{(m+1)}_u$ has at least $m\cdot (m+1)^{u-2}$ proper prefixes of the form $T_t^{(m)}$: there is one prefix for every combination of exponents $a_1,\dots,a_{u-1}\in\{0,\dots,m\}$ with $a_1\neq0$. Therefore $T_{m\cdot (m+1)^{u-2}}^{(m)}$ is a prefix of $B_u^{(m+1)}$. This proves the first inequality.
	
	Lemma~\ref{Ttm characterization}~\ref{a} tells us that the string $B_u^{(m+1)}$ has at most $(m+1)^{u-1}$ proper prefixes of the form $T_t^{(m)}$. Therefore the string $T_{(m+1)^{u-1}+1}^{(m)}$ cannot be a proper prefix of $B_u^{(m+1)}$. Now $B_u^{(m+1)}$ must be a prefix of $T_{(m+1)^{u-1}+1}^{(m)}$. This proves the second inequality.
	
	In case that $u=m+3$, consider the $(m+1)^{u-1}$ proper prefixes from the previous paragraph. One of them is $B_{m+2}^{(m+1)}$. This string has suffix $T_{m+2}^{(m+1)}$. From Lemma~\ref{t strings property} it follows that $T_{m+2}^{(m+1)}=(m+2)^{m+1}$. Now $T_{m+2}^{(m+1)}$ has curling number $m+1$, so $B_{m+2}^{(m+1)}$ has curling number at least $m+1$. It follows by Lemma~\ref{t strings property} that $B_{m+2}^{(m+1)}$ cannot be of the form $T_t^{(m)}$. There are then only $(m+1)^{u-1}-1$ proper prefixes left, so the second inequality is strict.
\end{proof}	
\end{lemma}

 From these estimates we deduce two propositions that we use in further sections:
\begin{proposition}\label{estimates 2}
	For all $m\geq1$ and $t\geq2$ we have $$\beta^{(m+1)}(\lceil \log_{m+1}(t)\rceil)\leq \tau^{(m)}(t)\leq \beta^{(m+1)}(\lceil \log_{m+1}(t)\rceil+2).$$
\begin{proof}
Since $\tau^{(m)}$ is an increasing function, we have $$\tau^{(m)}(t)\geq  \tau^{(m)}( (m+1)^{\lceil \log_{m+1}(t)\rceil-1}+1).$$
Suppose that $\lceil \log_{m+1}(t)\rceil=1$. Then $$\tau^{(m)}((m+1)^{\lceil \log_{m+1}(t)\rceil-1}+1)=\tau^{(m)}(2)=1=\beta^{(m+1)}(1)=\beta^{(m+1)}(\lceil \log_{m+1}(t)\rceil).$$ Now suppose that $\lceil \log_{m+1}(t)\rceil\geq2$. Then we apply Lemma~\ref{estimates} with $u=\lceil \log_{m+1}(t)\rceil$. It follows that $\tau^{(m)}((m+1)^{\lceil \log_{m+1}(t)\rceil-1}+1)\geq \beta^{(m+1)}(\lceil \log_{m+1}(t)\rceil)$. In both cases we obtain that $\tau^{(m)}(t)\geq\beta^{(m+1)}(\lceil \log_{m+1}(t)\rceil)$.

On the other hand we have $$\tau^{(m)}(t)\leq  \tau^{(m)}(m\cdot (m+1)^{\lceil \log_{m+1}(t)\rceil}).$$ By Lemma~\ref{estimates} this is at most $\beta^{(m+1)}(\lceil \log_{m+1}(t)\rceil+2)$. 
\end{proof}
\end{proposition}
\begin{proposition}\label{estimates 3}
	For all $m\geq1$ and $n\geq m+3$, we have
	$$m\cdot (m+1)^{t(m+1,n)-2}\leq t(m,n)\leq (m+1)^{t(m+1,n)-1}+1.$$
	The second inequality is strict when $n=m+3$.
\begin{proof}
We apply Lemma~\ref{estimates} to $u= t(m+1,n)$ with $n\geq m+3$ and $m\geq1$ to see that
$$\tau^{(m)}(m\cdot (m+1)^{t(m+1,n)-2})\leq \beta^{(m+1)}(t(m+1,n))\leq \tau^{(m)}((m+1)^{t(m+1,n)-1}+1).$$
When $n=m+3$, we have $u=t(m+1,m+3)=m+3$ by Example~\ref{t example}, so then the second inequality is strict. 
Now using Equation~\ref{tau tau eq}, we get
$$\tau^{(m)}(m\cdot (m+1)^{t(m+1,n)-2})\leq \tau^{(m)}(t(m,n))\leq \tau^{(m)}((m+1)^{t(m+1,n)-1}+1).$$
Since $\tau^{(m)}$ is a strictly increasing function, the same inequalities hold for its arguments.
\end{proof}
\end{proposition}
In this last proposition we connect $t(m,n)$ and $t(m+1,n)$. The exponential nature of this connection explains the tower of exponents in Theorem~\ref{full tower}.

\section{Bounds for \texorpdfstring{$\beta^{(m)}(u)$}{beta}}\label{the asymptotic behavior of beta}

In this section we give an upper bound for $\beta^{(m)}(u)$, which is the last building block for our main theorem. Van de Bult et al.~\cite[p.\ 19]{MR2268453} observed from calculations that there are real constants $\varepsilon_m$ for $m\geq 1$ such that for each $m$, we have $\beta^{(m)}(u)\approx\varepsilon_m\cdot (m+1)^{u-1}$ as $u\rightarrow\infty$. We prove that this is indeed the case.

Intuitively, this can be explained as follows. Since $B^{(m)}_{u+1}=(B^{(m)}_u)^{m+1}S^{(m)}_u$, we have $$\beta^{(m)}(u+1)=(m+1)\cdot \beta^{(m)}(u)+\sigma^{(m)}(u).$$ If $\sigma^{(m)}(u)$ is relatively small, then it follows that $\beta^{(m)}(u+1)$ is approximately equal to $(m+1)\cdot\beta^{(m)}(u)$. Iterating this we obtain that $\beta^{(m)}(u)\approx (m+1)^{u-1}$, and the contributions from the $\sigma^{(m)}(u)$ give us the extra factor $\varepsilon_m$.

As a first step we can use the above equality to deduce that $(m+1)\cdot\beta^{(m)}(u)\leq\beta^{(m)}(u+1)$, so $\beta^{(m)}(u)\leq\frac{\beta^{(m)}(u+1)}{m+1}$. With induction we obtain the following:
\begin{equation}
\label{fractionbetainequality}
	\beta^{(m)}(t)\leq\frac{\beta^{(m)}(u)}{(m+1)^{u-t}}\text{ for all }t,u\text{ with }1\leq t\leq u.
\end{equation}
Then we see that
\begin{equation}
\label{nondecreasing}
	\frac{\beta^{(m)}(t)}{(m+1)^{t-1}}\leq\frac{\beta^{(m)}(u)}{(m+1)^{u-1}}\text{ for all }t,u\text{ with }1\leq t\leq u.
\end{equation}
Hence, if we can prove that $\frac{\beta^{(m)}(u)}{(m+1)^{u-1}}$ is bounded as $u\rightarrow\infty$, then it follows that $\frac{\beta^{(m)}(u)}{(m+1)^{u-1}}$ converges for $u\rightarrow\infty$. The limit value is then equal to the constant $\varepsilon_m$.

We use the inequalities from the previous section to obtain an upper bound for $\sigma^{(m)}$: for all $m,t\geq 1$ we have
\begin{align*}
	\sigma^{(m)}(t)&\leq \tau^{(m)}(t+1)\\
	&\leq \beta^{(m+1)}(\lceil\log_{m+1}(t+1)\rceil+2)&\text{ (by Proposition~\ref{estimates 2})},
\end{align*}
hence
\begin{equation}\label{keyineq}
	\sigma^{(m)}(t)\leq\beta^{(m+1)}(\lceil\log_{m+1}(t+1)\rceil+2).
\end{equation}
We are now ready to prove the upper bound for $\beta^{(m)}(u)$.
\begin{lemma}\label{beta fraction bounded}
	Consider the numbers $M_7,M_8,\dots\in\mathbb R$ defined by $M_7:=4$ and $M_{k+1}:= M_k\cdot(1+9\cdot(k+1)^{1.59}/2^k)$ for $k\geq7$. Then 
	\begin{equation*}\frac{\beta^{(m)}(u)}{(m+1)^{u-1}}<M_u\text{ for all }m\geq 1,u\geq7.
	\end{equation*} 
	\begin{proof}
		We prove by induction on $u$ that the lemma holds for all $u\geq7$. First consider $u=7$. From the values of $\beta^{(m)}(t)$ that van de Bult et al. calculated~\cite[Table 6]{MR2268453}, we see that $\beta^{(m)}(7)<4\cdot(m+1)^6$ for $m=1,2,3,4,5$. For $m\geq 6$, we know from Lemma 4.1 by van de Bult et al.~\cite[p.\ 19]{MR2268453} that $\beta^{(m)}(7)=\frac{(m+1)^7-1}{m}<4\cdot (m+1)^6$. So the lemma holds for $u=7$.

		Now suppose that for a fixed $k\geq7$, we have $\beta^{(m)}(k)<M_k\cdot(m+1)^{k-1}$ for all $m\geq 1$. Then we have:
		\begin{align*}
			\frac{\beta^{(m)}(k+1)}{(m+1)^k}&=\frac{(m+1)\cdot\beta^{(m)}(k)}{(m+1)^{k}}+\frac{\sigma^{(m)}(k)}{(m+1)^{k}}\\
			&\leq \frac{\beta^{(m)}(k)}{(m+1)^{k-1}}+\frac{\beta^{(m+1)}(\lceil\log_{m+1}(k+1)\rceil+2)}{(m+1)^k}&&\text{ (by Inequality~\ref{keyineq})}\\
			&< M_k+\frac{\beta^{(m+1)}(\lceil\log_{m+1}(k+1)\rceil+2)}{(m+1)^k}&&\text{ (induction hypothesis)}\\
			&\leq M_k+\frac{\beta^{(m+1)}(k)}{(m+2)^{k-\lceil\log_{m+1}(k+1)\rceil-2}\cdot(m+1)^k}&&\text{ (by Inequality~\ref{fractionbetainequality}; see below)}\\
			&<M_k+\frac{M_k\cdot (m+2)^{\lceil\log_{m+1}(k+1)\rceil+1}}{(m+1)^k}&&\text{ (induction hypothesis)}\\
			&\leq M_k\cdot\left(1+\frac{(m+1)^{\log_2(3)\cdot (\lceil\log_{m+1}(k+1)\rceil+1)}}{(m+1)^k}\right)&&\text{ (see below)}\\
			&\leq  M_k\cdot (1+(m+1)^{\log_2(3)\cdot (\lceil\log_2(k+1)\rceil+1)-k})\\
			&\leq M_k\cdot(1+2^{\log_2(3)\cdot (\lceil\log_2(k+1)\rceil+1)-k})&&\text{ (see below)}\\
			&< M_k\cdot(1+9\cdot(k+1)^{1.59}/2^k)&&\text{ (see below)}\\
			&=M_{k+1}.
		\end{align*}
		Here we used that the following four inequalities hold for all $k\geq 7$ and $m\geq1$: 
		\begin{align*}
		\lceil \log_{m+1}(k+1)\rceil+2&\leq k\\
		\log_{m+1}(m+2)&\leq \log_2(3)\\
		\log_2(3)\cdot (\lceil\log_2(k+1)\rceil+1)-k&\leq 0\\
		\log_2(3)&<1.59.
		\end{align*}
		Proving these inequalities is straightforward.
	\end{proof}
\end{lemma}

 From the lemma and Inequality~\ref{nondecreasing}, we obtain that $\frac{\beta^{(m)}(u)}{(m+1)^{u-1}}<M$ for all $m,u\geq 1$, where we define
$M$ as $4\cdot\prod_{k=7}^{\infty}(1+9\cdot(k+1)^{1.59}/2^k)$. Our calculations give $M= 96.504\dots$, so $M<100$. Hence 
\begin{equation}\label{100}
	\frac{\beta^{(m)}(u)}{(m+1)^{u-1}}<100
\end{equation}
for all $m,u\geq 1$.

In the end we only need bounds for $\beta^{(1)}(u)$. A lower bound is easy: Inequality~\ref{fractionbetainequality} gives us that $\beta^{(1)}(u)\geq 2^{u-1}\cdot\beta^{(1)}(1)=2^{u-1}$. Hence
\begin{equation}2^{u-1}\leq \beta^{(1)}(u)<100\cdot 2^{u-1}.\label{beta1 bounds}
\end{equation}

\section{The tower of exponents}\label{the tower of exponents}
We combine Proposition~\ref{estimates 3} and Inequality~\ref{beta1 bounds} to prove our main result, which we restate for convenience:

\bigskip
\noindent{\bf Theorem 1.}
	For all integers $n\geq4$, the first occurrence of $n$ in Gijswijt's sequence is at position
\begin{equation*}2^{2^{\scriptstyle 3^{\scriptstyle 4^{\Ddots^{^{\scriptstyle (n-2)^{\scriptstyle \alpha}}}}}}},
\end{equation*}
		for a real number $\alpha\in (n-2,n-1)$ that depends on $n$.

	\begin{proof}
	The first occurrence of $n=4$ is at position $220$. This equals $2^{2^{\scriptstyle 2.96\dots}}$.
	
	Now suppose that $n\geq 5$.	The location of the first $n$ is at position $\beta^{(1)}(t(1,n))$. By Inequality~\ref{beta1 bounds} this is at least $2^{t(1,n)-1}$. We can now repeatedly apply the first part of Proposition~\ref{estimates 3}, which states that $$t(m,n)\geq m\cdot (m+1)^{t(m+1,n)-2}$$ for all $m\geq1$ and $n\geq m+3$. Using this for $m=1$ yields $2^{t(1,n)-1}\geq 2^{2^{\scriptstyle t(2,n)-2}-1}\geq 2^{2^{\scriptstyle t(2,n)-3}}$. We also see that $$t(m,n)-3\geq (m+1)^{t(m+1,n)-3}$$ for $2\leq m\leq n-4$. Combining this we obtain
	$$2^{2^{\scriptstyle t(2,n)-3}}\geq 	2^{2^{\scriptstyle 3^{\scriptstyle 4^{\Ddots^{^{\scriptstyle  (n-3)^{\scriptstyle t(n-3,n)-3}}}}}}}.$$
	For $m=n-3$, we have $$t(n-3,n)-3\geq (n-3)\cdot (n-2)^{t(n-2,n)-2}-3=(n-3)\cdot(n-2)^{n-2}-3>(n-2)^{n-2}.$$ 
	Here we used Example~\ref{t example} which implies that $t(n-2,n)=n$.
	We conclude that the position of the first $n$ is larger than 
	$$2^{2^{\scriptstyle 3^{	\Ddots^{^{\scriptstyle (n-3)^{\scriptstyle (n-2)^{\scriptstyle n-2}}}}}}}.$$
		
	On the other hand, Inequality~\ref{beta1 bounds} tells us that $\beta^{(1)}(t(1,n))<100\cdot2^{t(1,n)-1}$. We can now use the second part of Proposition~\ref{estimates 3}, which gives us that $$t(m,n)-1\leq (m+1)^{t(m+1,n)-1}$$
	for all $m\geq1$ and $n\geq m+3$. By repeatedly using this
	we see that $100\cdot 2^{t(1,n)-1}$ is at most $$100\cdot 2^{2^{\scriptstyle 3^{	\Ddots^{^{\scriptstyle (n-3)^{{\scriptstyle t(n-3,n)-1}}}}}}}.$$ We apply the proposition again for $m=n-3$, and this time it is a strict inequality.
	Therefore we get the upper bound $$100\cdot 2^{2^{\scriptstyle 3^{	\Ddots^{^{\scriptstyle (n-3)^{\left[\scriptstyle (n-2)^{\scriptstyle t(n-2,n)-1}-1\right]}}}}}}=100\cdot 2^{2^{\scriptstyle 3^{	\Ddots^{^{\scriptstyle (n-3)^{\left[\scriptstyle (n-2)^{\scriptstyle n-1}-1\right]}}}}}}.$$
	Since $n\geq5$ this is smaller than $$2^{2^{\scriptstyle 3^{	\Ddots^{^{\scriptstyle (n-2)^{\scriptstyle n-1}}}}}}.$$
	\end{proof}

There is some room for improvement:
\begin{remark}\label{narrow down}
	It is possible to narrow down the interval of $\alpha$, thereby giving tighter bounds for the position of the first $n$. This can be done by constructing a closed-form expression for $t(n-3,n)$. Then Proposition~\ref{estimates 3} needs to be used one less time for obtaining the tower of exponents.
\end{remark}
\begin{remark}\label{higher order}
We could also generalize Theorem~\ref{full tower} to the sequences $A^{(m)}$ for $m\geq2$. That would look something like this: For $n$ large enough, the first occurrence of $n$ in $A^{(m)}$ is at position $$(m+1)^{(m+1)^{\scriptstyle (m+2)^{\bigcdot\thinspace \ ^{\bigcdot\thinspace \ ^{\bigcdot^{\scriptstyle (n-2)^{\scriptstyle \alpha}}}}}}},$$
for a real number $\alpha\in(n-2,n-1)$ that depends on $m$ and $n$. We already have all the ingredients for this: both Proposition~\ref{estimates 3} and Inequality~\ref{100} hold for all $m$.
\end{remark}

\section{Acknowledgments}
The author would like to thank Gunther Cornelissen and Dion Gijswijt for supervising his bachelor thesis, which eventually led to this paper. He is grateful to the anonymous referees, whose insightful comments and suggestions significantly improved the manuscript. Finally, he thanks Brendan Miller for his helpful comments and for proofreading the English.

\bigskip
\hrule
\bigskip

\noindent 2020 {\it Mathematics Subject Classification}:
Primary 11B37; Secondary 11B83.

\noindent \emph{Keywords: } curling number, Gijswijt's sequence, recurrence, growth rate.

\bigskip
\hrule
\bigskip

\noindent (Concerned with sequences
\seqnum{A090822},
\seqnum{A091409},
\seqnum{A091411},
\seqnum{A091579},
\seqnum{A091787},
\seqnum{A091799},
\seqnum{A091840}, and
\seqnum{A357064}.)
\end{document}